\documentclass[12pt]{article}
\usepackage{amssymb,amsmath,amsthm,amsfonts,enumerate, bbm, mathdots}

\usepackage{latexsym}
\usepackage{amscd}
\usepackage{stmaryrd}
\usepackage{color}
\usepackage[all]{xypic}
\usepackage{epsfig}
\usepackage{graphics}
\usepackage{ifthen}
\usepackage{varioref}
\usepackage{rotating}
\usepackage{extarrows}
\usepackage{cite}
\usepackage{mathrsfs}

\numberwithin{equation}{section}

\theoremstyle{plain}
\newtheorem{thm}{Theorem}[section]
\newtheorem{cor}[thm]{Corollary}

\newtheorem{prop}[thm]{Proposition}

\newtheorem{eg}[thm]{Example}

\newmuskip\pFqmuskip
\newcommand*\pFq[6][8]{%
  \begingroup 
  \pFqmuskip=#1mu\relax
  \mathchardef\normalcomma=\mathcode`,
  \mathcode`\,=\string"8000
  \begingroup\lccode`\~=`\,
  \lowercase{\endgroup\let~}\pFqcomma
  {}_{#2}F_{#3}{\left(\genfrac..{0pt}{}{#4}{#5};#6\right)}%
  \endgroup
}
\newcommand{\pFqcomma}{{\normalcomma}\mskip\pFqmuskip}

\makeatletter
\newenvironment{proofof}[1]{\par
  \pushQED{\qed}%
  \normalfont \topsep6\p@\@plus6\p@\relax
  \trivlist
  \item[\hskip\labelsep
        \bfseries
    Proof of #1\@addpunct{.}]\ignorespaces
}{%
  \popQED\endtrivlist\@endpefalse
}
\makeatother

\definecolor{darkgreen}{rgb}{0.0625,0.64,0.0625}
\usepackage[
            pdfstartview=FitH,
            CJKbookmarks=true,
            bookmarksnumbered=true,
            bookmarksopen=true,
            colorlinks,
            pdfborder=001,
            linkcolor=blue,
            anchorcolor=green,
            citecolor=red
            ]{hyperref}

\newfont{\scyr}{wncyr10 scaled 550}

\def\proof{\noindent {\bf Proof.\;}}

\def\wt{\operatorname{wt}}
\def\dep{\operatorname{dep}}
\def\height{\operatorname{ht}}

\allowdisplaybreaks

\begin{document}

\title{Interpolated multiple $t$-values of general level with fixed weight, depth and height}
\date{~}

\author{{Zhonghua Li${}^{a,}$\thanks{Email: zhonghua\_li@tongji.edu.cn}\quad and \quad{Zhenlu Wang${}^{b,}$\thanks{Email: zhenluwang@tzc.edu.cn}}}\\[1mm]
\small a. School of Mathematical Sciences,\\ \small Key Laboratory of Intelligent Computing and Applications (Ministry of Education), \\ \small Tongji University, Shanghai 200092, China\\
\small b. School of Electronics and Information Engineering, \\ \small Taizhou University, Taizhou, Zhejiang 318000, China}

\maketitle
\begin{abstract}
In this paper, we introduce the interpolated multiple $t$-values of general level and represent a generating function for sums of interpolated multiple $t$-values of general level with fixed weight, depth, and height in terms of a generalized hypergeometric function $\,_3F_2$ evaluated at $1$. Furthermore, we explore several special cases of our results. The theorems presented in this paper extend earlier results on multiple zeta values and multiple $t$-values of general level.
\end{abstract}

{\small
{\bf Keywords}  multiple zeta values; multiple $t$-values; interpolated multiple $t$-values; generalized hypergeometric function.

{\bf 2020 Mathematics Subject Classification} 11M32, 33C20
}

\section{Introduction}

The theory of multiple zeta values (MZVs) originated in the early 1990s, independently developed by Hoffman~\cite{Hoffman92} and Zagier~\cite{Zagier1994}. These values, which Euler first studied in special cases, have since played a significant role in both number theory and mathematical physics. For an index $\mathbf{k} = (k_1, \ldots, k_n)$ with $k_i \in \mathbb{Z}_{\geq 1}$, we define three key quantities:
\begin{itemize}
    \item the weight: $\wt(\mathbf{k}) = k_1 + \cdots + k_n$,
    \item the depth: $\dep(\mathbf{k}) = n$,
    \item the height: $\height(\mathbf{k}) = |\{i \mid 1 \leq i \leq n, \, k_i \geq 2\}|$.
\end{itemize}
An index $\mathbf{k}$ is called admissible if $k_1 > 1$. For such an index, the MZV~\cite{Hoffman92,Zagier1994}  is defined by the following  convergent series
$$
\zeta(\mathbf{k}) = \zeta(k_1, \ldots, k_n) = \sum_{m_1 > \cdots > m_n > 0} \frac{1}{m_1^{k_1} \cdots
 m_n^{k_n}}.$$
A closely related variant, the multiple zeta-star values (MZSVs), has also been widely investigated. For an admissible index $\mathbf{k}$, the MZSV~\cite{Hoffman92}
 is given by 
$$ 
\zeta^\star(\mathbf{k}) = \zeta^\star(k_1, \ldots, k_n) = \sum_{m_1 \geq \cdots \geq m_n > 0} \frac{1}{m_1^{k_1} \cdots
 m_n^{k_n}}.
$$
It is well-known that numerous $\mathbb{Q}$-linear relations exist among MZVs and MZSVs; see~\cite{Zhao2016} for further details. Nevertheless, their exact algebraic structure remains largely mysterious, and many profound questions in this area await resolution.

Hoffman~\cite{Hoffman2019} introduced and studied an odd variant of MZVs (resp. MZSVs), namely the multiple $t$-values (MtVs) (resp. the multiple $t$-star values (MtSVs)). Recent research has revealed that MtVs and MtSVs possess rich arithmetic and algebraic structures remarkably similar to those of their classical counterparts. Moreover, these odd variants appear to be deeply connected to the algebraic structure of the space generated by all MZVs. For detailed investigations of these connections, we refer to~
\cite{Hoffman2019,Kaneko-Tsumura,Li-Wang,Murakami2021}.

We considered the  level $N$ variants in \cite{Li-Wang}. Let $N$ be a fixed positive integer. For an admissible index $\mathbf{k}=(k_1,\ldots,k_n)$ and any $a\in \{1,2,\ldots,N\}$, the MtV of level $N$ and the MtSV of level $N$ are defined by
$$t_{N,a}(\mathbf{k})=t_{N,a}(k_1,\ldots,k_n)=\sum\limits_{m_1>\cdots>m_n>0\atop m_i\equiv a\pmod N}\frac{1}{m_1^{k_1}\cdots m_n^{k_n}}$$
and
$$t_{N,a}^{\star}(\mathbf{k})=t_{N,a}^{\star}(k_1,\ldots,k_n)=\sum\limits_{m_1\geq\cdots\geq m_n>0\atop m_i\equiv a\pmod N}\frac{1}{m_1^{k_1}\cdots m_n^{k_n}},$$
respectively. Note that if the index is empty, we treat the values $t_{N,a}(\varnothing)=t_{N,a}^\star(\varnothing)=1$.
In the case of $N=2$ and $a=1$, we recover the classical MtVs~\cite{Hoffman2019} 
$$t(k_1,\ldots,k_n)=t_{2,1}(k_1,\ldots,k_n)=\sum\limits_{m_1>\cdots>m_n>0\atop m_i: \text{odd}}\frac{1}{m_1^{k_1}\cdots m_n^{k_n}}$$
and the MtSVs~\cite{Hoffman2019} 
$$t^\star(k_1,\ldots,k_n)=t^\star_{2,1}(k_1,\ldots,k_n)=\sum\limits_{m_1\geq\cdots\geq m_n>0\atop m_i: \text{odd}}\frac{1}{m_1^{k_1}\cdots m_n^{k_n}}$$
respectively. Moreover, it is easy to see that $t_{1,1}(\mathbf{k})=\zeta(\mathbf{k})$ and $t_{1,1}^{\star}(\mathbf{k})=\zeta^{\star}(\mathbf{k})$
which are MZVs and MZSVs respectively. We also have $t_{N,N}(\mathbf{k})=N^{-\wt({\mathbf{k}})}\zeta(\mathbf{k})$ and $t_{N,N}^\star(\mathbf{k})=N^{-\wt({\mathbf{k}})}\zeta^{\star}(\mathbf{k})$.

Yamamoto~\cite{Yamamoto} introduced the interpolated MZVs, which can be regarded as interpolation polynomials of MZVs and MZSVs. Li~\cite{Li2022} established a generating function for sums of interpolated MZVs with fixed weight, depth, and height, which generalizes Ohno-Zagier relation for MZVs~\cite{Ohno-Zagier} and the analogous result for MZSVs proved by Aoki, Kombu and Ohno~\cite{Aoki-Kombu-Ohno}. In this paper, we study interpolation polynomials of MtVs and MtSVs of general level. Let $r$ be a variable. For an admissible index $\mathbf{k}=(k_1,\ldots,k_n)$, the interpolated MtV of level $N$ is defined as
\begin{align*}
t_{N,a}^r(\mathbf{k})=t_{N,a}^r(k_1,\ldots,k_n)=\sum\limits_{\mathbf{p}}r^{n-\dep(\mathbf{p})}t_{N,a}(\mathbf{p}),
\end{align*}
where the summation extends over all indices $\mathbf{p}$ obtained from the pattern
$$(k_1 \Box k_2 \Box \cdots \Box  k_n)$$
by replacing each placeholder $\Box$ with either a comma ``$,$'' or a plus ``$+$ ''. Note that $t_{N,a}^0(\mathbf{k})=t_{N,a}(\mathbf{k})$ and $t_{N,a}^1(\mathbf{k})=t_{N,a}^\star(\mathbf{k})$. The parameter $r$ thus continuously interpolates between the non-star and star variants of MtVs of level $N$, providing a unified framework for studying these objects.

Our main object of study is the generating function for sums of interpolated MtVs of level $N$ with fixed weight, depth, and height. This generating function can be expressed in terms of the generalized hypergeometric function $_{3}F_2$ evaluated at $1$.  Here for a positive integer $m$ and complex numbers $b_1,\ldots,b_{m+1}$, $c_1,\ldots,c_m$ with none of $c_i$ is zero or a negative integer, the generalized hypergeometric function $_{m+1}F_m$ is defined by
\begin{align*}
\pFq{m+1}{m}{b_1,\ldots,b_{m+1}}{c_1,\ldots,c_m}{z}=\sum\limits_{n=0}^\infty\frac{(b_1)_n\cdots(b_{m+1})_n}{(c_1)_n\cdots(c_m)_n}\frac{z^n}{n!},
\end{align*}
where the Pochhammer symbol $(b)_n$ denotes the rising factorial:
$$(b)_n=\frac{\Gamma(b+n)}{\Gamma(b)}=\begin{cases}
1 & \text{if}\quad n=0, \\
b(b+1)\cdots(b+n-1) & \text{if}\quad n>0.
\end{cases}$$
The above series is absolutely and uniformly convergent for $|z|<1$ and the convergence also extends over the unit circle if $\Re\left(\sum c_i-\sum b_i\right)>0$.

This paper is organized as follows. We obtain a representation of the generating function for sums of interpolated MtVs of level $N$ with fixed weight, depth, and height in Section \ref{Sec:Main theorems}. Based on our main theorem, we discuss several special cases in Section \ref{Sec:Applications} and get a generating function for sums of interpolated MtVs of level $N$ with maximal height and a weighted sum formula of sums of interpolated MtVs of level $N$ with fixed weight and depth.

\section{Main theorems}\label{Sec:Main theorems}
Let $I_0(k,n,s)$ be the set of all admissible indices of weight $k$, depth $n$ and height $s$. Notice that $I_0(k,n,s)$ is nonempty if and only if $k\geq n+s$ and $n\geq s\geq1$. Let
\begin{align*}
X_0^r(k,n,s)=X_{N,a,0}^r(k,n,s)=\sum\limits_{\mathbf{k}\in I_0(k,n,s)}t_{N,a}^r(\mathbf{k})
\end{align*}
be the sum of interpolated MtVs of level $N$ with weight $k$, depth $n$, and height $s$. Then a generating function for these sums can be represented by the generalized hypergeometric function $_3F_2$ evaluated at $1$.

\begin{thm}\label{Thm:r-mtv generating function}
For formal variables $u,v,w$, we have
\begin{align*}
\sum\limits_{k\geq n+s\atop n\geq s\geq1}X_0^r(k,n,s)u^{k-n-s}v^{n-s}w^{2s-2}=\frac{1}{(a-\alpha_1)(a-\alpha_2)}\pFq{3}{2}{\frac{a+\beta_1}{N},\frac{a+\beta_2}{N},1}{\frac{a-\alpha_1}{N}+1,\frac{a-\alpha_2}{N}+1}{1},
\end{align*}
where $\alpha_1,\alpha_2$ are determined by $\alpha_1+\alpha_2=u+vr, \alpha_1\alpha_2=r(uv-w^2)$, and $\beta_1,\beta_2$ are determined by $\beta_1+\beta_2=-u+v(1-r), \beta_1\beta_2=(r-1)(uv-w^2)$.
\end{thm}

Theorem \ref{Thm:r-mtv generating function} unifies our previous results on MtVs and MtSVs of level $N$ and extends the corresponding result on interpolated MZVs. Specifically: 
\begin{itemize}
  \item[(i)] The case $r=0$ recovers \cite[Theorem 2.1]{Li-Wang} for the MtVs of level $N$; 
  \item[(ii)] The case $r=1$ recovers \cite[Theorem 2.2]{Li-Wang} for the MtSVs of level $N$;  
  \item[(iii)] The case $N=a=1$ reduces to \cite[Theorem 3.1]{Li2022} concerning the interpolated MZVs through application of the transformation formula \cite[7.4.4.1]{Prudnikov-Brychkov-Marichev}.
\end{itemize} 

We prove the theorem similarly as that in \cite[Theorems 2.1, 2.2]{Li-Wang}. For an index $\mathbf{k}=(k_1,\ldots,k_n)$, we define
\begin{align*}
\mathcal{L}_{N,a}(\mathbf{k};z)=\sum\limits_{m_1>\cdots>m_n>0\atop m_i\equiv a\pmod N}\frac{z^m_1}{m_1^{k_1}\cdots m_n^{k_n}},
\end{align*}
which converges absolutely for $|z|<1$. It is easy to see that $\mathcal{L}_{N,a}(\mathbf{k};1)=t_{N,a}(\mathbf{k})$ for any admissible index $\mathbf{k}$, which is a MtV of level $N$. Define that
\begin{align*}
\mathcal{L}_{N,a}^r(\mathbf{k};z)=\sum\limits_{\mathbf{p}}r^{n-\dep(\mathbf{p})}\mathcal{L}_{N,a}(\mathbf{p};z),
\end{align*}
where $\mathbf{p}$ runs over all sequences of the form
$$\mathbf{p}=(k_1\Box k_2\Box\cdots\Box k_n)$$
in which each $\Box$ is filled by the comma ``$,$'' or the plus ``$+$''. Then the right-hand side of the definition equation of $\mathcal{L}_{N,a}^r(\mathbf{k};z)$ converges locally uniformly in the domain $|z|<1$, and converges in points $|z|=1$ if $\mathbf{k}$ is admissible with the value of it coincides with $t_{N,a}^r(\mathbf{k})$. We have the following iterated integral representation
\begin{align*}
\mathcal{L}_{N,a}^r(\mathbf{k};z)=\int\limits_{z>z_1>\cdots>z_k>0}\prod\limits_{i=1}^kf_i(z_i)dz_i,
\end{align*}
where $k=\wt(\mathbf{k})=k_1+\cdots+k_n$ and
\begin{align*}
f_i(z)=\begin{cases}
\frac{r}{z}+\frac{z^{N-1}}{1-z^N} &  \text{if}\quad i=k_1,k_1+k_2,\ldots,k_1+\cdots+k_{n-1}, \\
\frac{z^{a-1}}{1-z^N} & \text{if}\quad i=k, \\
\frac{1}{z} & \text{otherwise}.
\end{cases}
\end{align*}
It is easy to find that
\begin{align}\label{Eq:derivation equation}
\frac{d}{dz}\mathcal{L}_{N,a}^r(k_1,\ldots,k_n;z)=\begin{cases}
\frac{1}{z}\mathcal{L}_{N,a}^r(k_1-1,k_2,\ldots,k_n;z) &  \text{if}\quad k_1>1, \\
\left(\frac{r}{z}+\frac{z^{N-1}}{1-z^N}\right)\mathcal{L}_{N,a}^r(k_2,\ldots,k_n;z) & \text{if}\quad k_1=1, n\geq 2, \\
\frac{z^{a-1}}{1-z^N} & \text{if}\quad k_1=n=1.
\end{cases}
\end{align}

For nonnegative integers $k,n,s$, we denote by $I(k,n,s)$ the set of all indices of weight $k$, depth $n$ and height $s$. We define the following two sums
\begin{align*}
X^r(k,n,s;z)=\sum\limits_{\mathbf{k}\in I(k,n,s)}\mathcal{L}_{N,a}^r(\mathbf{k};z)\quad\text{and}\quad X_0^r(k,n,s;z)=\sum\limits_{\mathbf{k}\in I_0(k,n,s)}\mathcal{L}_{N,a}^r(\mathbf{k};z).
\end{align*}
One can see that $X^r(k,n,s;z)$ and $X_0^r(k,n,s;z)$ also depend on $N$ and $a$, but we do not add extra subscripts $N,a$ for brevity. By convention, the above sums are defined to be zero if the index set is empty. Additionally, we set $X^r(0,0,0;z)=1$.

For integers $k,n,s$, by \eqref{Eq:derivation equation}, we obtain
\begin{itemize}
\item[(1)] if $k\geq n+s$ and $n\geq s\geq1$, then
\begin{align}\label{Eq:dX_0}
&\frac{d}{dz}X_0^r(k,n,s;z)=\frac{1}{z}\left[X^r(k-1,n,s-1;z)+X_0^r(k-1,n,s;z)\right.\notag\\
&\qquad\qquad\qquad\qquad\left.-X_0^r(k-1,n,s-1;z)\right];
\end{align}
\item[(2)] if $k\geq n+s$, $n\geq s\geq0$ and $n\geq2$, then
\begin{align}\label{Eq:d(X-X_0)}
\frac{d}{dz}\left[X^r(k,n,s;z)-X_0^r(k,n,s;z)\right]=\left(\frac{r}{z}+\frac{z^{N-1}}{1-z^N}\right)X^r(k-1,n-1,s;z).
\end{align}
\end{itemize}

Now we define the generating functions for $X^r(k,n,s;z)$ and $X^r_0(k,n,s;z)$ by
$$\Phi^r(z)=\Phi^r(u,v,w;z)=\sum\limits_{k,n,s\geq0}X^r(k,n,s;z)u^{k-n-s}v^{n-s}w^{2s}$$ and
$$\Phi_{0}^r(z)=\Phi^r_{0}(u,v,w;z)=\sum\limits_{k,n,s\geq0}X^r_{0}(k,n,s;z)u^{k-n-s}v^{n-s}w^{2s-2},$$ respectively. Notice that we do not add extra subscripts $N,a$ to $\Phi^r(z)$ and $\Phi^r_0(z)$ for brevity. Using \eqref{Eq:derivation equation}, \eqref{Eq:dX_0} and \eqref{Eq:d(X-X_0)}, we have
\begin{align*}
&\frac{d}{dz}\Phi_{0}^r(z)=\frac{1}{vz}\left(\Phi^r(z)-1-w^2\Phi^r_{0}(z)\right)+\frac{u}{z}\Phi_{0}^r(z),\\
&\frac{d}{dz}\left(\Phi^r(z)-w^2\Phi_{0}^r(z)\right)=v\left(\frac{r}{z}+\frac{z^{N-1}}{1-z^N}\right)\left(\Phi^r(z)-1\right)+\frac{vz^{a-1}}{1-z^N}.
\end{align*}
Eliminating $\Phi^r(z)$, we obtain the differential equation satisfied by $\Phi_0^r(z)$.

\begin{prop}
$\Phi_0^r=\Phi_0^r(z)$ satisfies the following differential equation
\begin{align}\label{Eq:Phi_0^r differential}
&z^2(1-z^N)(\Phi_0^r)''+z[(1-u)(1-z^N)-v\left(r+(1-r)z^N\right)](\Phi_0^r)'\notag\\
&\qquad\qquad\qquad\qquad+\left(r+(1-r)z^N\right)(uv-w^2)\Phi_0^r=z^a.
\end{align}
\end{prop}

The solution to the differential equation \eqref{Eq:Phi_0^r differential} is presented below.
\begin{thm}\label{Thm:Phi_0^r rep}
  For formal variables $u,v,w$, we have
  \begin{align*}
  \Phi^r_{0}(u,v,w;z)=\frac{z^a}{(a-\alpha_1)(a-\alpha_2)}\pFq{3}{2}{\frac{a+\beta_1}{N},\frac{a+\beta_2}{N},1}{\frac{a-\alpha_1}{N}+1,\frac{a-\alpha_2}{N}+1}{z^N},
  \end{align*}
  where $\alpha_1,\alpha_2$ are determined by $\alpha_1+\alpha_2=u+vr, \alpha_1\alpha_2=r(uv-w^2)$, and $\beta_1,\beta_2$ are determined by $\beta_1+\beta_2=-u+v(1-r), \beta_1\beta_2=(r-1)(uv-w^2)$.
  \end{thm}
  \proof
Assume that $\Phi_0^r(z)=\sum\limits_{n=1}^\infty p_nz^n$, then we have $p_1,\ldots,p_{a-1},p_{a+1},\ldots,p_N=0$,
\begin{align*}
p_a=\frac{1}{a(a-u-vr)+r(uv-w^2)}
\end{align*}
and
\begin{align*}
p_{n+N}=\frac{n(n-u+v-vr)-(1-r)(uv-w^2)}{(n+N)(n+N-u-vr)+r(uv-w^2)}p_n,\quad n\geq1.
\end{align*}
Hence for any $m\geq1$ that satisfies $m\not\equiv a\pmod N$, we have $p_m=0$. And for any $n\geq 1$, we have
\begin{align*}
p_{a+nN}&=\frac{(a+(n-1)N)(a+(n-1)N-u+v-vr)-(1-r)(uv-w^2)}{(a+nN)(a+nN-u-vr)+r(uv-w^2)}p_{a+(n-1)N}.
\end{align*}
Then we get
\begin{align*}
p_{a+nN}&=\frac{(a+(n-1)N+\beta_1)(a+(n-1)N+\beta_2)}{(a+nN-\alpha_1)(a+nN-\alpha_2)}p_{a+(n-1)N}\\
&=\frac{\left(\frac{a+\beta_1}{N}\right)_n\left(\frac{a+\beta_2}{N}\right)_n}{\left(\frac{a-\alpha_1}{N}+1\right)_n\left(\frac{a-\alpha_2}{N}+1\right)_n}\frac{1}{(a-\alpha_1)(a-\alpha_2)}.
\end{align*}
Hence, we finish the proof.
\qed

Theorem \ref{Thm:Phi_0^r rep} provides an explicit expression for the generating function of $X^r_0(k,n,s;z)$ by using generalized hypergeometric function $_3F_2$. Now we prove Theorem \ref{Thm:r-mtv generating function}.

\begin{proofof}{Theorem \ref{Thm:r-mtv generating function}}
Setting $z=1$ in Theorem \ref{Thm:Phi_0^r rep}, since $\mathcal{L}_{N,a}^r(\mathbf{k};1)=t_{N,a}^r(\mathbf{k})$ for admissible indices $\mathbf{k}$, we get Theorem \ref{Thm:r-mtv generating function}.
\end{proofof}

\section{Applications}\label{Sec:Applications}
\subsection{Sum of height one}
We denote by $\{k\}^n$ the sequence of $k$ repeated $n$ times. Setting $w=0$ in Theorem \ref{Thm:r-mtv generating function}, then we have $\alpha_1=u, \alpha_2=vr$, and $\beta_1=-u, \beta_2=v(1-r)$. Therefore, we get a representation of the generating function for interpolated MtVs of level $N$ with height one.

\begin{cor}
For formal variables $u,v$, we have
\begin{align*}
\sum\limits_{k\geq n+1,n\geq1}t_{N,a}^r(k-n+1,\{1\}^{n-1})u^{k-n-1}v^{n-1}=\frac{1}{(a-u)(a-vr)}\pFq{3}{2}{\frac{a-u}{N},\frac{a+v(1-r)}{N},1}{\frac{a+N-u}{N},\frac{a+N-vr}{N}}{1}.
\end{align*}
\end{cor}

\cite[Lemma 2.7]{Li-Wang} shows that the right-hand side of above formula is reducible to power sums. Hence for any integer $m\geq2$, we have
\begin{align*}
\sum\limits_{n=1}^\infty t_{N,a}^r(m,\{1\}^{n-1})v^{n-1}=\frac{1}{a^{m-1}(a-vr)}\pFq{m+1}{m}{1,\frac{a+v(1-r)}{N},\{\frac{a}{N}\}^{m-1}}{\frac{a+N-vr}{N},\{\frac{a+N}{N}\}^{m-1}}{1}.
\end{align*}

\subsection{Sum of maximal height}
Setting $v=0$ in Theorem \ref{Thm:r-mtv generating function}, then we get $\alpha_1+\alpha_2=u,\alpha_1\alpha_2=-rw^2$, and $\beta_1+\beta_2=-u,\beta_1\beta_2=(1-r)w^2$. By using the summation formula \cite[7.4.4.28]{Prudnikov-Brychkov-Marichev}
\begin{align}\label{Eq:3F2 maximal height}
\pFq{3}{2}{a,b,1}{c,2+a+b-c}{1}=\frac{1+a+b-c}{(1+a-c)(1+b-c)}\left(1-c+\frac{\Gamma(c)\Gamma(1+a+b-c)}{\Gamma(a)\Gamma(b)}\right)
\end{align}
with $a=\frac{a+\beta_1}{N}$, $b=\frac{a+\beta_2}{N}$ and $c=\frac{a-\alpha_1}{N}+1$, we find that
\begin{align*}
\frac{1}{(a-\alpha_1)(a-\alpha_2)}\pFq{3}{2}{\frac{a+\beta_1}{N},\frac{a+\beta_2}{N},1}{\frac{a-\alpha_1}{N}+1,\frac{a-\alpha_2}{N}+1}{1}
=\frac{-1}{w^2}+\frac{1}{w^2}\frac{\Gamma\left(\frac{a-\alpha_1}{N}\right)\Gamma\left(\frac{a-\alpha_2}{N}\right)}{\Gamma\left(\frac{a+\beta_1}{N}\right)\Gamma\left(\frac{a+\beta_2}{N}\right)}.
\end{align*}

Hence by using \cite[Lemma 2.8]{Li-Wang}, we get the following generating function for sums of interpolated MtVs of level $N$ with maximal height.
\begin{cor}
For formal variables $u,w$, we have
\begin{align}\label{Eq:maximal height}
1+\sum\limits_{k\geq 2n,n\geq1}X_0^r(k,n,n)u^{k-2n}w^{2n}=\exp\left\{\sum\limits_{n=2}^\infty\frac{t_{N,a}(n)}{n}(\alpha_1^n+\alpha_2^n-\gamma_1^n-\gamma_2^n)\right\},
\end{align}
where $\alpha_1,\alpha_2$ are determined by $\alpha_1+\alpha_2=u,\alpha_1\alpha_2=-rw^2$, and $\gamma_1,\gamma_2$ are determined by $\gamma_1+\gamma_2=u,\gamma_1\gamma_2=(1-r)w^2$.
\end{cor}

Furthermore, setting $u=0$ in \eqref{Eq:maximal height}, we obtain a generating function for $t^r_{N,a}(\{2\}^n)$ displayed as
\begin{align*}
1+\sum\limits_{n=1}^\infty t^r_{N,a}(\{2\}^n)w^{n}=\exp\left(\sum\limits_{n=1}^\infty\frac{\left(r^n-(r-1)^n\right)t_{N,a}(2n)}{n}w^{n}\right).
\end{align*}

\subsection{A weighted sum formula}
Let $I_0(k,n)$ be the set of admissible indices of weight $k$ and depth $n$. We have the following generating function for weight sums of interpolated MtVs of level $N$.

\begin{cor}
  For a formal variable $u$, we have
  \begin{align}\label{Eq:GF weighted sum formula}
&\sum\limits_{k>n\geq1}\left((1-2r)^{k-n-1}2^{n-1}\sum\limits_{\mathbf{k}\in I_0(k,n)}t_{2a,a}^r(\mathbf{k})\right)u^{k-2}\notag\\
  &=a^{-2}2^{\frac{2u}{a}}t(2)\exp\left(\sum\limits_{n=2}^\infty\frac{4(2^{n-1}-1)t(n)}{na^n(2^n-1)}u^n\right).
  \end{align}
\end{cor}
\proof
Setting $uv=w^2$ in Theorem \ref{Thm:r-mtv generating function}, then we have $\alpha_1+\alpha_2=u+vr,\beta_1+\beta_2=-u+v(1-r)$ and $\alpha_1\alpha_2=\beta_1\beta_2=0$. Hence we get
\begin{align*}
  \sum\limits_{k>n\geq1}\left(\sum\limits_{\mathbf{k}\in I_0(k,n)}t_{2a,a}^r(\mathbf{k})\right)u^{k-n-1}v^{n-1}=\frac{1}{a(a-u-vr)}\pFq{3}{2}{\frac{a-u+v(1-r)}{N},\frac{a}{N},1}{\frac{a-u-vr}{N}+1,\frac{a}{N}+1}{1}.
\end{align*}
We discuss the special case of $N=2a$. Let $u=(1-2r)x$ and $v=2x$, we get
\begin{align*}
  \sum\limits_{k>n\geq1}\left((1-2r)^{k-n-1}2^{n-1}\sum\limits_{\mathbf{k}\in I_0(k,n)}t_{2a,a}^r(\mathbf{k})\right)x^{k-2}=\frac{1}{a(a-x)}\pFq{3}{2}{\frac{a+x}{2a},\frac{1}{2},1}{\frac{3a-x}{2a},\frac{3}{2}}{1}.
\end{align*}
Using the summation formula \cite[7.4.4.21]{Prudnikov-Brychkov-Marichev}
\begin{align*}
&\pFq{3}{2}{a,b,c}{1+a-b,1+a-c}{1}\\
&=\frac{\sqrt{\pi}}{2^a}\frac{\Gamma(1+a-b)\Gamma(1+a-c)\Gamma\left(1+\frac{a}{2}-b-c\right)}{\Gamma\left(\frac{1+a}{2}\right)\Gamma\left(1+\frac{a}{2}-b\right)\Gamma\left(1+\frac{a}{2}-c\right)\Gamma(1+a-b-c)},\quad\Re(a-2b-2c)>-2,
\end{align*}
with $a=1$, $b=\frac{1}{2}$ and $c=\frac{a+x}{2a}$, we obtain
\begin{align*}
\frac{1}{a(a-x)}\pFq{3}{2}{\frac{a+x}{2a},\frac{1}{2},1}{\frac{3a-x}{2a},\frac{3}{2}}{1}&=\frac{1}{a(a-x)}\frac{\sqrt{\pi}}{2}\frac{\Gamma(\frac{3}{2})\Gamma(\frac{3}{2}-\frac{x}{2a})\Gamma(\frac{1}{2}-\frac{x}{2a})}{\Gamma(1)\Gamma(1)\Gamma(1-\frac{x}{2a})\Gamma(1-\frac{x}{2a})}\\&=\frac{\pi}{8a^2}\frac{\Gamma\left(\frac{1}{2}-\frac{x}{2a}\right)^2}{\Gamma\left(1-\frac{x}{2a}\right)^2}.
\end{align*}
Note that $\zeta(n)=(1-2^{-n})^{-1}t(n)$ and $t(2)=\frac{\pi^2}{8}$, by using the duplication formula
\begin{align*}
\Gamma\left(\frac{1}{2}-\frac{z}{2}\right)=\frac{\sqrt{\pi}2^z\Gamma(1-z)}{\Gamma\left(1-\frac{z}{2}\right)},
\end{align*}
and the expansion
\begin{align*}
\Gamma(1-z)=\exp\left(\gamma z+\sum\limits_{n=2}^\infty\frac{\zeta(n)}{n}z^n\right),
\end{align*}
we finish the proof of \eqref{Eq:GF weighted sum formula}.
\qed

Moreover, we get a weighted sum formula by expanding the right-hand side of \eqref{Eq:GF weighted sum formula}.
\begin{cor}\label{Cor:weighted sum formula}
For any integer $k\geq2$, we have
\begin{align*}
&\sum\limits_{n=1}^{k-1}\left((1-2r)^{k-n-1}2^{n-1}\sum\limits_{\mathbf{k}\in I_0(k,n)}t_{2a,a}^r(\mathbf{k})\right)\\
&=\sum\limits_{n+n_1+\cdots+n_m=k-2\atop n,m\geq0,n_1,\ldots,n_m\geq2}\frac{2^{n+2m}(2^{n_1-1}-1)\cdots(2^{n_m-1}-1)}{a^kn!m!n_1\cdots n_m(2^{n_1}-1)\cdots(2^{n_m}-1)}t(2)t(n_1)\cdots t(n_m)\log^n(2).
\end{align*}
\end{cor}

\begin{eg}
Considering weight $k=3$ or $4$, we get the following weighted formulas for MtVs and MtSVs derived from Corollary \ref{Cor:weighted sum formula}.
For $k=3$, we have
$$t(3)+2t(2,1)=-t^\star(3)+2t^\star(2,1)=2t(2)\log(2).$$
For $k=4$, we have
\begin{align*}
  t(4)+2(t(3,1)+t(2,2))+4t(2,1,1)&=t^\star(4)-2(t^\star(3,1)+t^\star(2,2))+4t^\star(2,1,1)\\&=\frac{2}{3}t^2(2)+2t(2)\log^2(2).
\end{align*}
\end{eg}


\section*{Acknowledgments}

The first author is supported by the Natural Science Foundation of Shanghai (Grant No. 24ZR1469000).


\end{document}